\documentclass[11pt]{amsart}
\usepackage{mathrsfs}
\usepackage{amssymb,latexsym}
\usepackage{pstcol,pstricks,color}

 \numberwithin{equation}{section}

\newtheorem{theorem}{Theorem}[section]
\newtheorem{proposition}[theorem]{Proposition}

\newtheorem{definition}[theorem]{Definition}

\newtheorem{remark}[theorem]{Remark}

\newtheorem{question}[theorem]{Question}

\title{Pattern Avoidance in Generalized Non-crossing Trees}

\begin{document}
\maketitle
\begin{center}
Yidong Sun\footnote{Corresponding author: Yidong Sun,
sydmath@yahoo.com.cn.}$^\dag$ and Zhiping Wang$^\ddag$

Department of Mathematics, Dalian Maritime University, 116026 Dalian, P.R. China\\[5pt]

{\it $^\dag$sydmath@yahoo.com.cn,  $^\ddag$wangzhiping5006@tom.com }
\end{center}\vskip0.5cm

\subsection*{Abstract} In this paper, the problem of pattern avoidance in
generalized non-crossing trees is studied. The generating functions
for generalized non-crossing trees avoiding patterns of length one
and two are obtained. Lagrange inversion formula is used to obtain
the explicit formulas for some special cases. Bijection is also
established between generalized non-crossing trees with special
pattern avoidance and the little Schr\"{o}der paths.

\medskip

{\bf Keywords}: Generalized non-crossing tree, pattern avoidance,
Catalan number, little Schr\"{o}der path.

\noindent {\sc 2000 Mathematics Subject Classification}: Primary
05A05; Secondary 05C30

\section{Introduction}

A {\em non-crossing tree } (NC-tree for short) is a tree drawn on
$n$ points in $\{1,2,\cdots,n\}$ numbered in counterclockwise order
on a circle such that the edges lie entirely within the circle and
do not cross. Non-crossing trees have been investigated by Chen and
Yan \cite{ChenYan}, Deutsch and Noy \cite{DeutNoy}, Flajolet and Noy
\cite{FlajNoy}, Gu, et al. \cite{GuLi}, Hough \cite{Hough}, Noy
\cite{Noy}, Panholzer and Prodinger \cite{PanhProd}. Recently, some
problems of pattern avoidance in NC-trees have been studied by Sun
and Wang \cite{SunWang}. It is well known that the set of NC-trees
with $n+1$ vertices is counted by the generalized Catalan number
$\frac{1}{2n+1}\binom{3n}{n}$ \cite[A001764]{Sloane}.

A {\em generalized non-crossing tree } (GNC-tree for short) is a
modified NC-tree such that the labels are weakly increasing in
counterclockwise order and if $j\geq 1$ is a label then all $1\leq
i\leq j$ are also labels. See Figure \ref{fDD2} for example.

In the sequel, we are concerned with the rooted GNC-trees such that
the first 1 is the root. Let $\textsf{GNC}_n $ denote the set of
rooted GNC-trees of $n+1$ vertices. It is easy to prove that
$\textsf{GNC}_n $ is counted by
$|\textsf{GNC}_n|=\frac{2^n}{2n+1}\binom{3n}{n}$ \cite[not
listed]{Sloane}, for a bar can be or not be inserted into any
position between $i$ and $i+1$ for $1\leq i\leq n$ in an NC-tree of
$n+1$ points and assume that one bar always appears in the position
between $n+1$ and $1$, then relabel the numbers between any two bars
in a proper way to form a GNC-tree.

A {\em descent (an\ ascent, a level)} is an edge $(i,j)$ such that
$i>j\ (i<j, i=j)$ and $i$ is on the path from the root to the vertex
$j$. If encoding an ascent by $u$, a level by $h$ and a descent by
$d$, then each path in a GNC-tree can be represented by a ternary
word on $\{u,h,d\}$ by viewing from the root. In analogy with the
well-established permutation patterns \cite{SimSch, Wilf1}, we
propose a definition of patterns in GNC-trees.

\begin{definition}
Let $w=w_1w_2\dots w_n$ and $\sigma=\sigma_1\sigma_2\dots\sigma_k$
be two ternary words on $\{u,h,d\}$. Then $w$ contains the pattern
$\sigma$ if it has a subword $w_{i+1}w_{i+2}\dots w_{i+k}$ equal to
$\sigma$ for some $0\leq i\leq n-k$; otherwise $w$ is called
$\sigma$-avoiding. A {\rm{GNC}}-tree $T$ is called $\sigma$-avoiding
if $T$ has no subpath (viewing from the root) encoded by $\sigma$.
\end{definition}

Let $\mathcal{P}_k$ denote the set of ternary words of length $k$ on
$\{u,h,d\}$. For any $\sigma\in \mathcal{P}_k$, let
$\textsf{GNC}_n^m(\sigma)$ denote the set of GNC-trees in
$\textsf{GNC}_n$ which contain the pattern $\sigma$ exactly $m$
times. For any nonempty subset $P\subset\mathcal{P}_k$,
$\textsf{GNC}_n(P)$ denotes the set of GNC-trees in $\textsf{GNC}_n$
which avoid all the patterns in $P$. Analogous to restricted
permutations, a counterpart in GNC-trees is the following question
\begin{question}
Determine the cardinalities of $\textsf{GNC}_n(P)$ for $P\subset
\mathcal {P}_{k}$ and $\textsf{GNC}_n^m(\sigma)$ for $\sigma\in
\mathcal {P}_{k}$.
\end{question}

In the literature, two kind of special NC-trees have been
considered, that is non-crossing increasing trees and non-crossing
alternating trees. Both of them are counted by the Catalan numbers
$C_n=\frac{1}{n+1}\binom{2n}{n}$ \cite[A000108]{Sloane}. {\em A
non-crossing increasing (alternating) tree} is an NC-tree with the
vertices on the path from the root 1 to any other vertex appearing
in increasing (alternating) order. By our notation, a non-crossing
increasing tree is just a $d$-avoiding NC-tree and a non-crossing
alternating tree is just a $\{uu,dd\}$-avoiding NC-tree. Bijections
between non-crossing alternating trees and Dyck paths have been
presented in \cite{Stanley}. But for GNC-trees, it seems to be
thrown little light on this subject.

In this paper, we deal with several patterns and find the
corresponding generating functions for GNC-trees. More precisely, we
investigate the patterns in $\mathcal {P}_1$ in Section 2 and the
patterns in $\mathcal {P}_2$ in Section 3. Lagrange inversion
formula is used to obtain the explicit formulas for some special
cases. Bijection is also established between GNC-trees with special
pattern avoidance and the little Schr\"{o}der paths.

\section{The patterns in $\mathcal{P}_1$}

For any $T\in \textsf{GNC}_n$, let $u(T), h(T), d(T)$ denote the
number of ascents, levels and descents of $T$ respectively, then
$u(T)+h(T)+d(T)=n$. Let $\textsf{GNC}_n^{*}$ be the set of GNC-trees
$T$ in $\textsf{GNC}_n$ such that $T$ has only one point with label
1, namely only the root has the label 1 and others have labels
greater than 1. Define
\begin{eqnarray*}
T_{x,y,z}(t)&=&\sum_{n\geq 0}t^n\sum_{T\in
\textsf{GNC}_n}x^{u(T)}y^{h(T)}z^{d(T)} ,\\
T_{x,y,z}^{*}(t)&=&\sum_{n\geq 0}t^n\sum_{T\in
\textsf{GNC}_n^{*}}x^{u(T)}y^{h(T)}z^{d(T)}.
\end{eqnarray*}

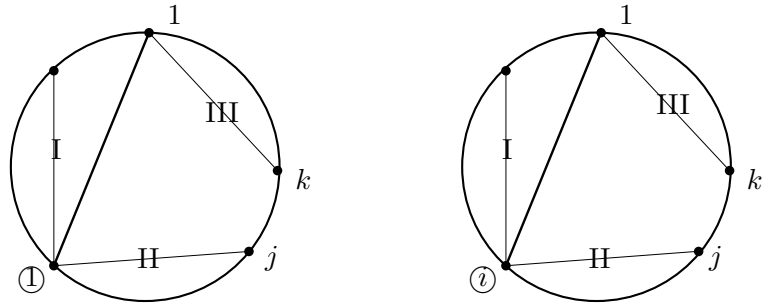
\begin{figure}[h]
\setlength{\unitlength}{0.5mm}
\begin{pspicture}(0,0)(3,3)
\pscircle(-1,1){1.8} \pscircle(5,1){1.8}

\pscircle*[unit=18pt,linewidth=.5pt](-1.5,4.4){0.1}
\pscircle*[unit=18pt,linewidth=.5pt](-3.5,-.5){0.1}
\pscircle*[unit=18pt,linewidth=.5pt](0.6,-.2){0.1}
\pscircle*[unit=18pt,linewidth=.5pt](1.2,1.5){0.1}
\pscircle*[unit=18pt,linewidth=.5pt](-3.5,3.6){0.1}

\psline[unit=18pt,linewidth=.9pt](-1.5,4.4)(-3.5,-.5)
\psline*[unit=18pt,linewidth=.5pt](-3.5,-.5)(0.6,-.2)
\psline*[unit=18pt,linewidth=.5pt](-3.5,-.5)(-3.5,3.6)
\psline*[unit=18pt,linewidth=.5pt](-1.5,4.4)(1.2,1.5)

\put(-.7,2.9){$1$}\put(-2.7,-.6){\textcircled{$1$}}\put(0.6,-.3){$j$}\put(1,0.7){$k$}
\put(-0.2,1.6){III}\put(-2.25,1.1){I}\put(-1.1,-.35){II}

\pscircle*[unit=18pt,linewidth=.5pt](8,4.4){0.1}
\pscircle*[unit=18pt,linewidth=.5pt](6,-.5){0.1}
\pscircle*[unit=18pt,linewidth=.5pt](10.05,-.2){0.1}
\pscircle*[unit=18pt,linewidth=.5pt](10.7,1.5){0.1}
\pscircle*[unit=18pt,linewidth=.5pt](6,3.6){0.1}

\psline[unit=18pt,linewidth=.9pt](8,4.4)(6,-.5)
\psline*[unit=18pt,linewidth=.5pt](6,-.5)(10,-.2)
\psline*[unit=18pt,linewidth=.5pt](6,-.5)(6,3.6)
\psline*[unit=18pt,linewidth=.5pt](8,4.4)(10.7,1.5)

\put(5.3,2.9){$\tiny{1}$}\put(3.3,-.6){\textcircled{$i$}}
\put(6.5,-.3){$j$}\put(7,0.7){$k$}
\put(5.8,1.7){III}\put(3.75,1.1){I}\put(4.9,-.35){II}

\end{pspicture}\vskip0.5cm

\caption{ Decomposition of GNC-trees.}\label{fDD1}
\end{figure}

Close relations between $T_{x,y,z}(t)$ and $T_{x,y,z}^{*}(t)$ can be
established according to the decomposition of GNC-trees in Figure
\ref{fDD1}. Find the first and minimal label $i$, denoted by
\textcircled{$i$}, of $T\in \textsf{GNC}_n$ in counterclockwise
order such that the root $1$ and \textcircled{$i$} form an edge,
then $T$ can be partitioned into three parts.

\begin{itemize}
\item[(i)] The case $i=j=1$ and $k=1\ \mbox{or}\ 2$. Part I and II both
avoid the patterns $u$ and $d$, Part III still forms a GNC-tree.
Then the edge $(1,1)$ contributes an $yt$, each of Part I and II
contributes $T_{0,y,0}(t)$, and Part III contributes $T_{x,y,z}(t)$.

\item[(ii)] The case $i=1,j\geq 2$ and $k=j\ \mbox{or}\ j+1$. Part I
avoids the patterns $u$ and $d$, Part II still forms a GNC-tree
different from Part I. Then the edge $(1,1)$ contributes an $yt$,
Part I contributes $T_{0,y,0}(t)$ and Part II contributes
$T_{x,y,z}(t)-T_{0,y,0}(t)$. For Part III, except for the one point
case, decreasing all the labels (excluding the root), by $j-2$ units
(if $k=j$) or $j-1$ units (if $k=j+1$), one can obtain two GNC-trees
in $\textsf{GNC}_n^{*}$ for some $n\geq 1$. Then Part III
contributes $2{T}_{x,y,z}^{*}(t)-1$.

\item[(iii)] The case $i\geq 2, j\geq i$ and $k=j\ \mbox{or}\ j+1$.
The edge (1,\textcircled{$i$}) contributes an $xt$. Note that the
ascents and descents are exchanged in Part I, and the labels are
lying in $\{1,2,\dots, i\}$ or $\{2,3,\dots,i\}$, so Part I
contributes $2T_{z,y,x}(t)-T_{0,y,0}(t)$ for $i\geq 2$. Part II
still forms a GNC-tree after reducing the labels $\{i,\dots,j\}$ to
$\{1,\cdots,j-i+1\}$, so Part II contributes $T_{x,y,z}(t)$. Similar
to {\mbox(ii)}, Part III contributes $2T_{x,y,z}^{*}(t)-1$.

\end{itemize}
Summarizing these, we have
\begin{eqnarray}\label{eqn 2.1}
T_{x,y,z}(t)&=&1+ytT_{0,y,0}(t)^2T_{x,y,z}(t) \nonumber\\
            & &+\ ytT_{0,y,0}(t)\Big\{T_{x,y,z}(t)-T_{0,y,0}(t)\Big\}\Big\{2T_{x,y,z}^{*}(t)-1\Big\} \\
            & &+\ xtT_{x,y,z}(t)\Big\{2T_{z,y,x}(t)-T_{0,y,0}(t)\Big\}\Big\{2T_{x,y,z}^{*}(t)-1\Big\}.\nonumber
\end{eqnarray}
For any $T\in\textsf{GNC}_n^{*}$, by the similar decomposition, one
can derive that
\begin{eqnarray}\label{eqn 2.2}
T_{x,y,z}^{*}(t)=1+xtT_{x,y,z}(t)T_{z,y,x}(t)\Big\{2T_{x,y,z}^{*}(t)-1\Big\}.
\end{eqnarray}
Solve (\ref{eqn 2.2}) for $T_{x,y,z}^{*}(t)$ and substitute it into
(\ref{eqn 2.1}), one can get
\begin{eqnarray}\label{eqn 2.3}
T_{x,y,z}(t)&=&1-ytT_{0,y,0}(t)^2-xtT_{x,y,z}(t)T_{0,y,0}(t)+ytT_{0,y,0}(t)\Big\{1+T_{0,y,0}(t)\Big\}T_{x,y,z}(t) \nonumber \\
            & &+\
            2xt\Big\{1-ytT_{0,y,0}(t)^2\Big\}T_{x,y,z}(t)^2T_{z,y,x}(t).
\end{eqnarray}
Let $x=z=0$, (\ref{eqn 2.3}) reduces to
\begin{eqnarray}\label{eqn 2.4}
T_{0,y,0}(t)&=&1+ytT_{0,y,0}(t)^3.
\end{eqnarray}

Multiplying by $T_{0,y,0}(t)$ in both side of (\ref{eqn 2.3}) and
using (\ref{eqn 2.4}), after some routine computations, one can
deduce that
\begin{eqnarray}\label{eqn 2.5}
T_{x,y,z}(t)&=&1+(y-x)tT_{0,y,0}(t)^2T_{x,y,z}(t)+2xtT_{x,y,z}(t)^2T_{z,y,x}(t).
\end{eqnarray}
Solve (\ref{eqn 2.5}) for $T_{x,y,z}(t)$, we have
\begin{eqnarray}\label{eqn 2.51}
T_{x,y,z}(t)&=&\frac{1-(y-x)tT_{0,y,0}(t)^2-\sqrt{(1-(y-x)tT_{0,y,0}(t)^2)^2-8xtT_{z,y,x}(t)}}{4xtT_{z,y,x}(t)} \nonumber\\
            &=&\frac{1}{1-(y-x)tT_{0,y,0}(t)^2}C(\frac{2xtT_{z,y,x}(t)}{(1-(y-x)tT_{0,y,0}(t)^2)^2}),
\end{eqnarray}
where $C(t)=\frac{1-\sqrt{1-4t}}{2t}$ is the generating function for
Catalan numbers.

Exchanging $x$ and $z$ in (\ref{eqn 2.51}), one can obtain
$T_{z,y,x}(t)$, and then substitute it into (\ref{eqn 2.51}), one
has the following proposition.
\begin{proposition}
The generating functions $T_{x,y,z}(t)$ for GNC-trees satisfies
\begin{eqnarray*}
T_{x,y,z}(t)&=&\alpha C\big(2xt\alpha^2\beta
C\big(2zt\beta^2T_{x,y,z}(t)\big)\big),
\end{eqnarray*}
where $\alpha=\frac{1}{1-(y-x)tT_{0,y,0}(t)^2}$ and
$\beta=\frac{1}{1-(y-z)tT_{0,y,0}(t)^2}$.
\end{proposition}

By Lagrange inversion formula and some series expansions, one can
obtain the coefficients of $t^n$ of $T_{x,y,z}(t)$, but it seems to
be somewhat complicated. Now we will consider several special cases
which lead to interesting results.

\subsection{$u$-avoiding GNC-trees} Note that $T_{0,y,z}(t)$ is the
generating function for GNC-trees with no ascent. Let $x=0$ in
(\ref{eqn 2.5}), one can find
\begin{eqnarray*}
T_{0,y,z}(t)&=&1+ytT_{0,y,0}(t)^2T_{0,y,z}(t),
\end{eqnarray*}
from which, together with (\ref{eqn 2.4}), using Lagrange inversion
formula, one can deduce that
\begin{eqnarray*}
T_{0,y,z}(t)=T_{0,y,0}(t)=\sum_{n\geq
0}\frac{1}{2n+1}\binom{3n}{n}(yt)^n=T_{0,1,0}(yt).
\end{eqnarray*}
In fact, the above relation can be easily derived from the
definition of GNC-tree, for a $u$-avoiding GNC-tree must also avoid
the pattern $d$, such GNC-trees can be obtained by changing each
label of the underlying NC-trees to the label $1$.

\subsection{$h$-avoiding GNC-trees} Note that $T_{x,0,z}(t)$ is the
generating function for GNC-trees with no level. Let $y=0$ in
(\ref{eqn 2.4}) and (\ref{eqn 2.5}), one can find
\begin{eqnarray*}
T_{x,0,z}(t)&=&1-xtT_{x,0,z}(t)+2xtT_{x,0,z}(t)^2T_{z,0,x}(t),
\end{eqnarray*}
which, when $x=z=1$, generates
\begin{eqnarray*}
T_{1,0,1}(t)&=&1-tT_{1,0,1}(t)+2tT_{1,0,1}(t)^3.
\end{eqnarray*}
Let $\lambda=T_{1,0,1}(t)-1$, then
$\lambda=t(1+\lambda)(2(1+\lambda)^2-1)$, using Lagrange inversion
formula \cite{Wilf2}, one can deduce for $n\geq 1$ that
\begin{eqnarray*}
[t^n]T_{1,0,1}(t)&=&[t^n]\lambda=\frac{1}{n}[\lambda^{n-1}](1+\lambda)^n(2(1+\lambda)^2-1)^n
\\
&=&\frac{1}{n}\sum_{i=0}^n(-1)^{n-i}\binom{n}{i}2^i[\lambda^{n-1}](1+\lambda)^{n+2i}
\\
&=&\sum_{i=0}^n(-1)^{n-i}\frac{2^i}{2i+1}\binom{n}{i}\binom{n+2i}{n}
\\
&=&\sum_{i=0}^n(-1)^{n-i}\frac{2^i}{2i+1}\binom{3i}{i}\binom{n+2i}{3i}.
\end{eqnarray*}
Hence we have
\begin{theorem}
The set of $h$-avoiding GNC-trees of $n+1$ points is counted by
\begin{eqnarray*}
|\textsf{GNC}_n(h)|&=&
\sum_{i=0}^n(-1)^{n-i}\frac{2^i}{2i+1}\binom{3i}{i}\binom{n+2i}{3i}.
\end{eqnarray*}
\end{theorem}
This sequence beginning with $1,1,5,31,217,1637,12985$ is not listed
in Sloane's \cite{Sloane}.

\subsection{$d$-avoiding GNC-trees} Note that $T_{x,y,0}(t)$ is the
generating function for GNC-trees with no descent. Let $z=0$ in
(\ref{eqn 2.51}), we can find
\begin{eqnarray}
T_{x,y,0}(t)&=&
\frac{1}{1-(y-x)tT_{0,y,0}(t)^2}C(\frac{2xtT_{0,y,x}(t)}{(1-(y-x)tT_{0,y,0}(t)^2)^2}),\nonumber \\
&=&
\frac{1}{1-(y-x)tT_{0,1,0}(yt)^2}C(\frac{2xtT_{0,1,0}(yt)}{(1-(y-x)tT_{0,1,0}(yt)^2)^2}),
\label{eqn 2.6}
\end{eqnarray}
where we use the relation $T_{0,y,x}(t)=T_{0,y,0}(t)=T_{0,1,0}(yt)$.

\subsection*{Case i} Let $x=y=1$ in (\ref{eqn 2.6}), one gets
\begin{eqnarray*}
T_{1,1,0}(t)&=&C(2tT_{0,1,0}(t)).
\end{eqnarray*}
Taking the coefficient $t^n$ of $T_{1,1,0}(t)$, one gets
\begin{eqnarray*}
[t^n]T_{1,1,0}(t)&=&[t^n]C(2tT_{0,1,0}(t))=[t^n]\sum_{i\geq
0}2^iC_it^iT_{0,1,0}(t)^{i} \\
&=&[t^n]\sum_{i\geq 0}2^iC_it^i\sum_{j\geq
0}\frac{i}{3j+i}\binom{3j+i}{j}t^j,\\
&=&\sum_{i+j=n}\frac{i}{3j+i}\binom{3j+i}{j}2^iC_i.
\end{eqnarray*}
Hence we have
\begin{theorem}
The set of $d$-avoiding GNC-trees of $n+1$ points is counted by
\begin{eqnarray*}
|\textsf{GNC}_n(d)|&=&
\sum_{i+j=n}\frac{i}{3j+i}\binom{3j+i}{j}2^iC_i.
\end{eqnarray*}
\end{theorem}
This sequence beginning with $1,2,10,62,424,3070$ is not listed in
Sloane's \cite{Sloane}.

\subsection*{Case ii} Let $y=1$ in (\ref{eqn 2.6}), by (\ref{eqn 2.4}), one gets
\begin{eqnarray*}
T_{x,1,0}(t)&=&\frac{T_{0,1,0}(t)}{1+xtT_{0,1,0}(t)^3}C(\frac{2xtT_{0,1,0}(t)^3}{(1+xtT_{0,1,0}(t)^3)^2}).
\end{eqnarray*}
Taking the coefficient $t^nx^k$ of $T_{x,1,0}(t)$, one gets
\begin{eqnarray*}
[t^nx^k]T_{x,1,0}(t)&=&[t^{n-k}x^k]\frac{T_{0,1,0}(t)}{1+xT_{0,1,0}(t)^3}C(\frac{2xT_{0,1,0}(t)^3}{(1+xT_{0,1,0}(t)^3)^2})\\
&=&[t^{n-k}]\sum_{i=0}^{k} 2^iC_iT_{0,1,0}(t)^{3i+1}[x^{k-i}]\frac{1}{(1+xT_{0,1,0}(t))^{2i+1}} \\
&=&\sum_{i=0}^{k} 2^iC_i(-1)^{k-i}\binom{k+i}{k-i}[t^{n-k}]T_{0,1,0}(t)^{k+2i+1}  \\
&=&\sum_{i=0}^{k}
(-1)^{k-i}\binom{k+i}{k-i}\frac{k+2i+1}{3n-2k+2i+1}\binom{3n-2k+2i+1}{n-k}2^iC_i.
\end{eqnarray*}
Hence we have
\begin{theorem}
The set of $d$-avoiding GNC-trees of $n+1$ points with $k$ ascents
is counted by
\begin{eqnarray*}
\sum_{i=0}^{k}
(-1)^{k-i}\binom{k+i}{k-i}\frac{k+2i+1}{3n-2k+2i+1}\binom{3n-2k+2i+1}{n-k}2^iC_i.
\end{eqnarray*}
\end{theorem}

\subsection*{Case iii} 
\begin{figure}[h]\vskip0.2cm
\setlength{\unitlength}{0.5mm}
\begin{pspicture}(0,0)(3,3)
\pscircle(2,1){1.8}

\pscircle*[unit=18pt,linewidth=.5pt](3,4.4){0.1}
\pscircle*[unit=18pt,linewidth=.5pt](1.25,-.5){0.1}
\pscircle*[unit=18pt,linewidth=.5pt](5.3,-.2){0.1}
\pscircle*[unit=18pt,linewidth=.5pt](6,1.5){0.1}
\pscircle*[unit=18pt,linewidth=.5pt](1.25,3.6){0.1}
\pscircle*[unit=18pt,linewidth=.5pt](.35,1.4){0.1}
\pscircle*[unit=18pt,linewidth=.5pt](5.1,3.6){0.1}
\pscircle*[unit=18pt,linewidth=.5pt](3.3,-1.2){0.1}

\psline[unit=18pt,linewidth=.5pt](3,4.4)(1.25,-.5)
\psline[unit=18pt,linewidth=.5pt](1.25,-.5)(5.3,-.2)
\psline[unit=18pt,linewidth=.5pt](1.25,-.5)(3.3,-1.2)
\psline[unit=18pt,linewidth=.5pt](3,4.4)(1.25,3.6)
\psline[unit=18pt,linewidth=.5pt](3,4.4)(0.4,1.4)
\psline[unit=18pt,linewidth=.5pt](3,4.4)(6,1.5)
\psline[unit=18pt,linewidth=.5pt](6,1.5)(5.1,3.6)

\put(2.,2.9){1}\put(.5,-.5){$2$} \put(.5,2.3){2}\put(-.1,.9){$2$}
\put(2.1,-1.2){3}\put(3.5,-.5){$3$} \put(4,.9){4}\put(3.4,2.3){$5$}

\end{pspicture}\vskip1cm

\begin{pspicture}(10,3.5)
\psset{xunit=20pt,yunit=20pt}\psgrid[subgriddiv=1,griddots=10,
gridlabels=4pt](0,0)(16,3)

\psline(0,0)(1,1)(5,1)(6,2)(8,2)(10,0)(12,2)(14,0)

\pscircle*(0,0){0.06}\pscircle*(1,1){0.06}\pscircle*(3,1){0.06}\pscircle*(5,1){0.06}\pscircle*(6,2){0.06}
\pscircle*(8,2){0.06}\pscircle*(9,1){0.06}\pscircle*(10,0){0.06}\pscircle*(11,1){0.06}
\pscircle*(12,2){0.06}\pscircle*(13,1){0.06}\pscircle*(14,0){0.06}

\put(1.1,.9){\tiny$H_{\ell}H_r$}\put(2.5,.9){\tiny$H_{\ell}H_r$}\put(4.6,1.6){\tiny$H_{\ell}H_r$}

\put(5.5,2.6){$\Updownarrow$}

\end{pspicture}

\caption{ The bijection between $\textsf{GNC}_n(h,d)$ and
$\mathscr{R}_n$.}\label{fDD2}
\end{figure}
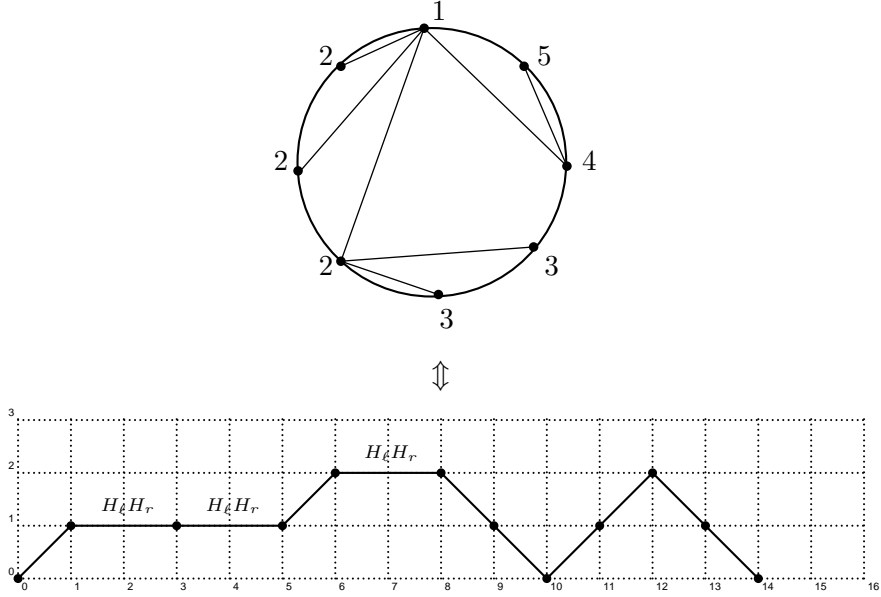

Let $x=1,y=0$ in (\ref{eqn 2.6}), using $T_{0,0,0}(t)=T_{0,1,0}(0)=1$, one gets
\begin{eqnarray*}
T_{1,0,0}(t)&=&\frac{1}{1+t}C(\frac{2t}{(1+t)^2})=\frac{1+t-\sqrt{1-6t+t^2}}{4t},
\end{eqnarray*}
which is the generating function for little Schr\"{o}der paths. A
{\em little Schr\"{o}der path of length $2n$} is a lattice path in
the first quadrant going from $(0,0)$ to $(2n,0)$ consisting of up
steps $U=(1,1)$, down steps $D=(1,-1)$ and horizontal steps
$H_{\ell}H_r=(2,0)$ with no horizontal step at the $x$-axis. Let
$\mathscr{R}_n$ denote the set of little Schr\"{o}der paths of
length $2n$ which is counted by the $n$th little Schr\"{o}der number
$R_n$ \cite[A001003]{Sloane}, whose generating function is
$R(t)=\frac{1+t-\sqrt{1-6t+t^2}}{4t}$. Hence we have
\begin{theorem}
The set of $\{h,d\}$-avoiding (i.e. increasing) GNC-trees of $n+1$
points is counted by the $n$th little Schr\"{o}der numbers. In other
words, there exists a bijection between $\textsf{GNC}_n(h,d)$ and
$\mathscr{R}_n$.
\end{theorem}
\begin{proof}
Read any $\{h,d\}$-avoiding GNC-tree of $n+1$ points in preorder,
denote an ascent $(i,j)$ by $U$ if it is read in the first time and
not following another ascent $(i,j)$; Denote an ascent $(i,j)$ by
$H_{\ell}$ if it is read in the second time and followed by another
ascent $(i,j)$ which is then denoted by $H_r$; Denote an ascent
$(i,j)$ by $D$ if it is read in the second time and not followed by
another ascent $(i,j)$. Then we can get a little Schr\"{o}der of
length $2n$. The above procedure is clearly invertible, see Figure
\ref{fDD2}.
\end{proof}

\section{The patterns in $\mathcal {P}_2$}
Let $\textsf{GNC}_n^{*}(\sigma)$ be the set of $\sigma$-avoiding
GNC-trees $T$ in $\textsf{GNC}_n^{*}$ with $\sigma\in \mathcal
{P}_2$. Define
\begin{eqnarray*}
T_{x,y,z}^{\sigma}(t)&=&\sum_{n\geq 0}t^n\sum_{T\in
\textsf{GNC}_n(\sigma)}x^{u(T)}y^{h(T)}z^{d(T)} ,\\
T_{x,y,z}^{*{\sigma}}(t)&=&\sum_{n\geq 0}t^n\sum_{T\in
\textsf{GNC}_n^{*}(\sigma)}x^{u(T)}y^{h(T)}z^{d(T)}.
\end{eqnarray*}
In this section, we will deal with the patterns $uu,dd,ud$ and $du$,
the others can be investigated similarly.

\subsection{The patterns $uu$ and $dd$.}
Close relations between $T_{x,y,z}^{uu}(t)$ and
$T_{x,y,z}^{*{uu}}(t)$ can be established according to the
decomposition of GNC-trees in Figure \ref{fDD1}.
\begin{itemize}

\item[(i)] The case $i=j=1$ and $k=1\ \mbox{or}\ 2$. Part I and II both
avoid the patterns $u$ and $d$, Part III still forms a $uu$-avoiding
GNC-tree. Then the edge $(1,1)$ contributes an $yt$, each of Part I
and II contributes $T_{0,y,0}(t)$, and Part III contributes
$T_{x,y,z}^{uu}(t)$.

\item[(ii)] The case $i=1,j\geq 2$ and $k=j\ \mbox{or}\ j+1$. Part I
avoids the patterns $u$ and $d$, Part II still forms a $uu$-avoiding
GNC-tree different from Part I. Then the edge $(1,1)$ contributes an
$yt$, Part I contributes $T_{0,y,0}(t)$ and Part II contributes
$T_{x,y,z}^{uu}(t)-T_{0,y,0}(t)$. For Part III, except for the one
point case, decreasing all the labels (excluding the root), by $j-2$
units if $k=j$ or $j-1$ units if $k=j+1$, one can obtain two
$uu$-avoiding GNC-trees in $\textsf{GNC}_n^{*}(uu)$ for some $n\geq
1$. Then Part III contributes $2{T}_{x,y,z}^{*uu}(t)-1$.

\item[(iii)] The case $i\geq 2, j\geq i$ and $k=j\ \mbox{or}\ j+1$.
The edge (1,\textcircled{$i$}) contributes an $xt$. Note that the
ascents and descents are exchanged in Part I, and the labels are
lying in $\{1,2,\dots, i\}$ or $\{2,3,\dots,i\}$, so Part I
contributes $2T_{z,y,x}^{dd}(t)-T_{0,y,0}(t)$. But Part II can not
begin with a $u$ edge, i.e., all edges (if exist) starting form
\textcircled{$i$} are $h$ edges, so we should further partition Part
II into three parts, by finding the last $i$, denoted by $i^{*}$, in
counterclockwise order such that (\textcircled{$i$},$i^{*}$) is an
$h$ edge which contributes an $yt$, see Figure \ref{fDD3}. Clearly,
Part $\mbox{II}_1$ and $\mbox{II}_2$ both contribute $T_{0,y,0}(t)$,
and Part $\mbox{II}_3$ contributes $T_{x,y,z}^{uu}(t)$. Similar to
{\mbox(ii)}, Part III contributes $2T_{x,y,z}^{*{uu}}(t)-1$.
\end{itemize}
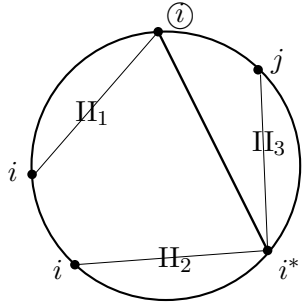
\begin{figure}[h]\vskip0.2cm
\setlength{\unitlength}{0.5mm}
\begin{pspicture}(0,0)(3,3)
\pscircle(2,1){1.8}

\pscircle*[unit=18pt,linewidth=.5pt](3,4.4){0.1}
\pscircle*[unit=18pt,linewidth=.5pt](1.25,-.5){0.1}
\pscircle*[unit=18pt,linewidth=.5pt](5.3,-.2){0.1}

\pscircle*[unit=18pt,linewidth=.5pt](.35,1.4){0.1}
\pscircle*[unit=18pt,linewidth=.5pt](5.1,3.6){0.1}

\psline[unit=18pt,linewidth=.9pt](3,4.4)(5.3,-.2)
\psline*[unit=18pt,linewidth=.5pt](1.25,-.5)(5.3,-.2)
\psline*[unit=18pt,linewidth=.5pt](5.3,-.2)(5.15,3.5)
\psline*[unit=18pt,linewidth=.5pt](3,4.4)(0.4,1.4)

\put(0.8,1.6){$\mbox{II}_1$}\put(1.9,-.35){$\mbox{II}_2$}\put(3.15,1.2){$\mbox{II}_3$}

\put(2.,2.9){\textcircled{$i$}}
\put(-.1,0.8){$i$}\put(.5,-.5){$i$}\put(3.5,-.5){$i^{*}$}
\put(3.4,2.3){$j$}

\end{pspicture}\vskip0.5cm

\caption{ The decomposition of Part II in (iii).}\label{fDD3}
\end{figure}

Summarizing these, we have
\begin{eqnarray}\label{eqn 3.1}
T_{x,y,z}^{uu}(t)&=&1+ytT_{0,y,0}(t)^2T_{x,y,z}^{uu}(t)+\ ytT_{0,y,0}(t)
\Big\{T_{x,y,z}^{uu}(t)-T_{0,y,0}(t)\Big\}\Big\{2T_{x,y,z}^{*{uu}}(t)-1\Big\}\nonumber \\
            & &+\ xt\Big\{2T_{z,y,x}^{dd}(t)-T_{0,y,0}(t)\Big\}\Big\{1+ytT_{0,y,0}(t)^2
            T_{x,y,z}^{uu}(t)\Big\}\Big\{2T_{x,y,z}^{*{uu}}(t)-1\Big\}.
\end{eqnarray}
For any $T\in\textsf{GNC}_n^{*}(uu)$, by the similar decomposition,
one can easily derive that
\begin{eqnarray}\label{eqn 3.2}
T_{x,y,z}^{*{uu}}(t)=1+xtT_{z,y,x}^{dd}(t)\Big\{1+ytT_{0,y,0}(t)^2T_{x,y,z}^{uu}(t)\Big\}\Big\{2T_{x,y,z}^{*}(t)-1\Big\}.
\end{eqnarray}
Similarly, for the pattern $dd$, close relations between
$T_{x,y,z}^{dd}(t)$ and $T_{x,y,z}^{*{dd}}(t)$ can be established
according to the decomposition of GNC-trees in Figure \ref{fDD1},
the details are omitted.
\begin{eqnarray}
T_{x,y,z}^{dd}(t)&=&1+ytT_{0,y,0}(t)^2T_{x,y,z}^{dd}(t)
            +ytT_{0,y,0}(t)\Big\{T_{x,y,z}^{dd}(t)-T_{0,y,0}(t)\Big\}\Big\{2T_{x,y,z}^{*{dd}}(t)-1\Big\} \nonumber\\
            & &+\
            xt\Big\{2T_{z,y,x}^{uu}(t)-T_{0,y,0}(t)\Big\}T_{x,y,z}^{dd}(t)
            \Big\{2T_{x,y,z}^{*{dd}}(t)-1\Big\}, \label{eqn 3.3}\\
\hskip0.4cm
T_{x,y,z}^{*{dd}}(t)&=&1+xtT_{z,y,x}^{uu}(t)T_{x,y,z}^{dd}(t)\Big\{2T_{x,y,z}^{*{dd}}(t)-1\Big\}.
\label{eqn 3.4}
\end{eqnarray}

Solve (\ref{eqn 3.2}) for $T_{x,y,z}^{*{uu}}(t)$ and (\ref{eqn 3.4})
for $T_{x,y,z}^{*{dd}}(t)$, and then substitute them respectively
into (\ref{eqn 3.1}) and (\ref{eqn 3.3}), after some
simplifications, one can get
\begin{eqnarray}
T_{x,y,z}^{uu}(t)&=&\Big\{1-xtT_{0,y,0}(t)^2+2xtT_{x,y,z}^{uu}(t)T_{z,y,x}^{dd}(t)\Big\}
                     \Big\{1+ytT_{0,y,0}(t)^2T_{x,y,z}^{uu}(t)\Big\},
        \label{eqn 3.5}\\
T_{x,y,z}^{dd}(t)&=&1+(y-x)tT_{0,y,0}(t)^2T_{x,y,z}^{dd}(t)
            +2xtT_{x,y,z}^{dd}(t)^2T_{z,y,x}^{uu}(t). \label{eqn 3.6}
\end{eqnarray}
Let $x=y=z=1$ in (\ref{eqn 3.5}) and (\ref{eqn 3.6}), we have
\begin{proposition}
The generating functions for $uu$-avoiding and $dd$-avoiding
GNC-trees are determined respectively by
\begin{eqnarray*}
T_{1,1,1}^{uu}(t)&=&\Big\{1-tT_{0,1,0}(t)^2+2tT_{1,1,1}^{uu}(t)T_{1,1,1}^{dd}(t)\Big\}
                     \Big\{1+tT_{0,1,0}(t)^2T_{1,1,1}^{uu}(t)\Big\}, \\
T_{1,1,1}^{dd}(t)&=&1+2tT_{1,1,1}^{dd}(t)^2T_{1,1,1}^{uu}(t).
\end{eqnarray*}
\end{proposition}

When $x=z=1,y=0$, by $T_{0,0,0}(t)=1$, (\ref{eqn 3.5}) and (\ref{eqn
3.6}) generate
\begin{eqnarray*}
T_{1,0,1}^{uu}(t)&=&1-t+2tT_{1,0,1}^{uu}(t)T_{1,0,1}^{dd}(t), \\
T_{1,0,1}^{dd}(t)&=&1-tT_{1,0,1}^{dd}(t)
+2tT_{1,0,1}^{dd}(t)^2T_{1,0,1}^{uu}(t),
\end{eqnarray*}
from which, one can deduce that
\begin{eqnarray}
T_{1,0,1}^{uu}(t)&=&\frac{1-t}{1-2tT_{1,0,1}^{dd}(t)}, \label{eqn 3.7}\\
T_{1,0,1}^{dd}(t)&=&1-3tT_{1,0,1}^{dd}(t)+4tT_{1,0,1}^{dd}(t)^2.
\label{eqn 3.8}
\end{eqnarray}
From (\ref{eqn 3.7}) and (\ref{eqn 3.8}), one can get
\begin{eqnarray*}
T_{1,0,1}^{dd}(t)&=&\frac{1+3t-\sqrt{(1+3t)^2-16t}}{8t} \\
&=& \frac{1}{1+3t}C(\frac{4t}{(1+3t)^2})=\sum_{i\geq 0}\frac{4^iC_it^i}{(1+3t)^{2i+1}} \\
&=& \sum_{i\geq 0}4^iC_it^i\sum_{j\geq 0}(-1)^j\binom{2i+j}{j}3^jt^j
\\
&=&\sum_{n\geq
0}t^n\sum_{j=0}^n(-1)^j\binom{2n-j}{j}3^j4^{n-j}C_{n-j}. \\
T_{1,0,1}^{uu}(t)&=&\frac{1-t}{1-2tT_{1,0,1}^{dd}(t)}=\frac{3-3t-\sqrt{(1+3t)^2-16t}}{2}\\
&=&\frac{3(1-t)-\sqrt{(1-t)^2-8t(1-t)}}{2} \\
&=&1-t+2tC(\frac{2t}{1-t})=1+t+\sum_{i\geq
0}\frac{2^{i+2}C_{i+1}t^{i+2}}{(1-t)^{i+1}} \\
&=&1+t+\sum_{n\geq 0}t^{n+2}\sum_{i=0}^n\binom{n}{i}2^{i+2}C_{i+1}.
\end{eqnarray*}
Hence we obtain
\begin{theorem}
The sets $\textsf{GNC}_{n+2}(uu,h)$ of $\{uu,h\}$-avoiding GNC-trees
and $\textsf{GNC}_{n}(dd,h)$ of $\{dd,h\}$-avoiding GNC-trees are
counted respectively by
\begin{eqnarray*}
|\textsf{GNC}_{n+2}(uu,h)|
&=&\sum_{i=0}^n\binom{n}{i}2^{i+2}C_{i+1},
\hskip3.2cm (\mbox{\cite[not\ listed]{Sloane}}),\\
|\textsf{GNC}_{n}(dd,h)|
&=&\sum_{j=0}^n(-1)^j\binom{2n-j}{j}3^j4^{n-j}C_{n-j}, \hskip1cm
(\mbox{\cite[A059231]{Sloane}}).
\end{eqnarray*}
\end{theorem}

\begin{remark}
Coker \cite{Coker} proved that $T_{1,0,1}^{dd}(t)$ is also the
generating function for $\mathscr{D}_n$, the set of different
lattice paths running from $(0,0)$ to $(2n,0)$ using steps from
$S=\{(k,\pm k): k\ \mbox{positive\ integer}\}$ that never go below
x-axis, and provided several different expressions for
$|\mathscr{D}_n|$. One can be asked to find a bijection between
$\mathscr{D}_n$ and $\textsf{GNC}_n(dd,h)$.
\end{remark}

\subsection{The patterns $ud$ and $du$} Similar to Subsection 3.1,
close relations between $T_{x,y,z}^{ud}(t)$ and
$T_{x,y,z}^{*{ud}}(t)$, and between $T_{x,y,z}^{du}(t)$ and
$T_{x,y,z}^{*{du}}(t)$, can be derived according to the
decomposition of GNC-trees, see in Figure \ref{fDD1}, but the
details are omitted.
\begin{eqnarray}
T_{x,y,z}^{ud}(t)&=&1+ytT_{0,y,0}(t)^2T_{x,y,z}^{ud}(t)
           + ytT_{0,y,0}(t)\Big\{T_{x,y,z}^{ud}(t)-T_{0,y,0}(t)\Big\}\Big\{2T_{x,y,z}^{*{ud}}(t)-1\Big\} \nonumber \\
            & &+\
            xtT_{x,y,z}^{ud}(t)\Big\{1+ytT_{0,y,0}(t)^2\big\{2T_{z,y,x}^{du}(t)-T_{0,y,0}(t)\big\}\Big\}
            \Big\{2T_{x,y,z}^{*{ud}}(t)-1\Big\}, \label{eqn 3.3.1}\\
\hskip0.6cm
T_{x,y,z}^{*{ud}}(t)&=&1+xtT_{x,y,z}^{ud}(t)\Big\{1+ytT_{0,y,0}(t)^2T_{z,y,x}^{du}(t)\Big\}\Big\{2T_{x,y,z}^{*{ud}}(t)-1\Big\},
\label{eqn 3.3.2} \\
T_{x,y,z}^{du}(t)&=&1+ytT_{0,y,0}(t)^2T_{x,y,z}^{du}(t)
            + ytT_{0,y,0}(t)\Big\{T_{x,y,z}^{du}(t)-T_{0,y,0}(t)\Big\}\Big\{2T_{x,y,z}^{*{du}}(t)-1\Big\} \nonumber\\
            & &+\
            xtT_{x,y,z}^{du}(t)\Big\{2T_{z,y,x}^{ud}(t)-T_{0,y,0}(t)\Big\}
            \Big\{2T_{x,y,z}^{*{du}}(t)-1\Big\},\label{eqn 3.3.3}\\
\hskip0.5cm
T_{x,y,z}^{*{du}}(t)&=&1+xtT_{x,y,z}^{du}(t)T_{z,y,x}^{ud}(t)\Big\{2T_{x,y,z}^{*{du}}(t)-1\Big\}.
\label{eqn 3.3.4}
\end{eqnarray}
Solve (\ref{eqn 3.3.2}) for $T_{x,y,z}^{*{ud}}(t)$ and (\ref{eqn
3.3.4}) for $T_{x,y,z}^{*{du}}(t)$, and then substitute them
respectively into (\ref{eqn 3.3.1}) and (\ref{eqn 3.3.3}), after
some simplifications, one can get
\begin{eqnarray}
\hskip0.5cm
T_{x,y,z}^{ud}(t)&=&1+(y-x)tT_{0,y,0}(t)^2T_{x,y,z}^{ud}(t)+
            2xtT_{x,y,z}^{ud}(t)^2\Big\{1+ytT_{0,y,0}(t)^2T_{z,y,x}^{du}(t)\Big\}, \label{eqn 3.3.5}
            \\
T_{x,y,z}^{du}(t)&=&1+(y-x)tT_{0,y,0}(t)^2T_{x,y,z}^{du}(t)+
            2xtT_{x,y,z}^{du}(t)^2T_{z,y,x}^{ud}(t). \label{eqn 3.3.6}
\end{eqnarray}
Let $x=y=z=1$ in (\ref{eqn 3.3.5}) and (\ref{eqn 3.3.6}), we have
\begin{proposition}
The generating functions for $ud$-avoiding and $du$-avoiding
GNC-trees are given by
\begin{eqnarray*}
T_{1,1,1}^{ud}(t)&=&1+2tT_{1,1,1}^{ud}(t)\Big\{1+tT_{0,1,0}(t)^2T_{1,1,1}^{du}(t)\Big\},           \\
T_{1,1,1}^{du}(t)&=&1+2tT_{1,1,1}^{du}(t)^2T_{1,1,1}^{ud}(t).
\end{eqnarray*}
\end{proposition}

When $x=z=1,y=0$, by $T_{0,0,0}(t)=1$, (\ref{eqn 3.3.5}) and
(\ref{eqn 3.3.6}) generate
\begin{eqnarray*}
T_{1,0,1}^{ud}(t)&=&1-tT_{1,0,1}^{ud}(t)+2tT_{1,0,1}^{ud}(t)^2, \\
T_{1,0,1}^{du}(t)&=&1-tT_{1,0,1}^{du}(t)
+2tT_{1,0,1}^{du}(t)^2T_{1,0,1}^{ud}(t),
\end{eqnarray*}
from which, one can deduce that
\begin{eqnarray*}
T_{1,0,1}^{ud}(t)&=&R(t)=\frac{1+t-\sqrt{1-6t+t^2}}{4t}=\frac{1}{1+t}C(\frac{2t}{(1+t)^2}), \\
T_{1,0,1}^{du}(t)&=&\frac{1+t-\sqrt{(1+t)^2-8tR(t)}}{4tR(t)} \\
&=&\frac{1}{1+t}C(\frac{2tR(t)}{(1+t)^2}) = \sum_{i\geq
0}\frac{2^iC_iR(t)^it^i}{(1+t)^{2i+1}} \\
&=&\sum_{i\geq 0}\frac{2^iC_it^i}{(1+t)^{3i+1}}\Big\{C(\frac{2t}{(1+t)^2})\Big\}^i \\
&=&\sum_{i\geq 0}\frac{2^iC_it^i}{(1+t)^{3i+1}}\sum_{j\geq
0}\frac{i}{2j+i}\binom{2j+i}{j}\frac{2^jt^j}{(1+t)^{2j}} \\
&=&\sum_{n\geq
0}t^n\sum_{i+j+k=n}(-1)^k\binom{3i+2j+k}{k}\frac{i}{2j+i}\binom{2j+i}{j}2^{i+j}C_i.
\end{eqnarray*}
Hence we have
\begin{theorem}
The set $\textsf{GNC}_n(ud,h)$ of $\{ud,h\}$-avoiding GNC-trees is
counted by the $n$th little Schr\"{o}der number, and the set
$\textsf{GNC}_{n}(du,h)$ of $\{du,h\}$-avoiding GNC-trees is counted
by
\begin{eqnarray*}
|\textsf{GNC}_{n}(du,h)| &=&
\sum_{i+j+k=n}(-1)^k\binom{3i+2j+k}{k}\frac{i}{2j+i}\binom{2j+i}{j}2^{i+j}C_i.
\end{eqnarray*}
\end{theorem}
This sequence beginning with $1,1,5,27,157,957,6025$ is not listed
in \cite{Sloane}.

\subsection{The pattern $\{uu, dd\}$} Now we consider the pattern $\{uu, dd\}$,
let $P=\{uu, dd\}$, according to the decomposition of GNC-trees,
relations between $T_{x,y,z}^{P}(t)$ and $T_{x,y,z}^{*P}(t)$ can be
derived, the details are omitted.
\begin{eqnarray*}
T_{x,y,z}^{P}(t)&=& 1+ ytT_{0,y,0}(t)^2T_{x,y,z}^{P}(t)+
ytT_{0,y,0}(t)\Big\{T_{x,y,z}^{P}(t)-T_{0,y,0}(t)\Big\} \Big\{2T_{x,y,z}^{*P}(t)-1\Big\}      \\
&&+\
xt\Big\{2T_{z,y,x}^{P}(t)-T_{0,y,0}(t)\Big\}\Big\{1+ytT_{0,y,0}(t)^2T_{x,y,z}^{P}(t)\Big\}
\Big\{2T_{x,y,z}^{*P}(t)-1\Big\} ,\\
T_{x,y,z}^{*P}(t)&=&1+xtT_{z,y,x}^{P}(t)\Big\{1+ytT_{0,y,0}(t)^2T_{x,y,z}^{P}(t)\Big\}\Big\{2T_{x,y,z}^{*P}(t)-1\Big\},
\end{eqnarray*}
from which, one can get
\begin{eqnarray}\label{eqn 3.4.0}
T_{x,y,z}^{P}(t)&=&
\Big\{1+ytT_{0,y,0}(t)^2T_{x,y,z}^{P}(t)\Big\}\Big\{1-xtT_{0,y,0}(t)^2+2xtT_{z,y,x}^{P}(t)T_{x,y,z}^{P}(t)\Big\}.
\end{eqnarray}
Let $x=y=z=1$ in (\ref{eqn 3.4.0}), one has
\begin{proposition}
The generating function for $\{uu,dd\}$-avoiding GNC-trees is given
by
\begin{eqnarray*}
T_{1,1,1}^{P}(t)&=&\Big\{1+tT_{0,1,0}(t)^2T_{1,1,1}^{P}(t)\Big\}
\Big\{1-tT_{0,1,0}(t)^2+2tT_{1,1,1}^{P}(t)^2\Big\}.
\end{eqnarray*}
\end{proposition}

Let $y=0$ in (\ref{eqn 3.4.0}), by $T_{0,0,0}(t)=1$, one can get
\begin{eqnarray}\label{eqn 3.4.1}
T_{x,0,z}^{P}(t)&=& 1-xt+2xtT_{z,0,x}^{P}(t)T_{x,0,z}^{P}(t).
\end{eqnarray}
Exchanging $x$ and $z$ in (\ref{eqn 3.4.1}), one has
\begin{eqnarray}\label{eqn 3.4.2}
T_{z,0,x}^{P}(t)&=& 1-zt+2ztT_{x,0,z}^{P}(t)T_{z,0,x}^{P}(t).
\end{eqnarray}
From (\ref{eqn 3.4.1}) and (\ref{eqn 3.4.2}), one can obtain
\begin{eqnarray*}
T_{x,0,z}^{P}(t)&=& 1-xt-2(z-x)tT_{x,0,z}^{P}(t)
+2ztT_{x,0,z}^{P}(t)^2,
\end{eqnarray*}
which leads to
\begin{eqnarray}\label{eqn 3.4.3}
T_{x,0,z}^{P}(t)&=&
\frac{1+2(z-x)t-\sqrt{(1+2(z-x)t)^2-8zt(1-xt)}}{4zt}.
\end{eqnarray}
Setting $z=1$ in (\ref{eqn 3.4.3}), one can deduce
\begin{eqnarray*}
T_{x,0,1}^{P}(t)&=&
\frac{1+2(1-x)t-\sqrt{(1+2(1-x)t)^2-8t(1-xt)}}{4t} \\
&=&\frac{1-xt}{1+2(1-x)t}C(\frac{2t(1-xt)}{(1+2(1-x)t)^2})=
\sum_{i\geq 0}\frac{2^iC_it^i(1-xt)^{i+1}}{(1+2(1-x)t)^{2i+1}} \\
&=&\sum_{i\geq
0}2^iC_it^i\sum_{j=0}^{i+1}(-1)^j\binom{i+1}{j}x^jt^j\sum_{k\geq
0}(-1)^k\binom{2i+k}{k}2^kt^k\sum_{\ell=0}^k(-1)^{\ell}\binom{k}{\ell}x^{\ell} \\
&=&\sum_{n\geq
0}\sum_{r=0}^nt^nx^r\sum_{i+j+k=n}(-1)^{r+k}\binom{i+1}{j}\binom{2i+k}{k}\binom{k}{r-j}2^{i+k}C_i.
\end{eqnarray*}
Setting $x=z=1$ in (\ref{eqn 3.4.3}), one has
\begin{eqnarray*}
T_{1,0,1}^{P}(t)&=& \frac{1-\sqrt{1-8t(1-t)}}{4t}=(1-t)C(2t(1-t)) \\
&=&\sum_{n\geq 0}t^n\sum_{i=0}^{n}(-1)^{n-i}\binom{i+1}{n-i}2^iC_i. \\
\end{eqnarray*}
Hence we have
\begin{theorem}
The number of $\{uu,dd,h\}$-avoiding (i.e., alternating) GNC-trees
of $n+1$ points is given by
\begin{eqnarray*}
\sum_{i=0}^{n}(-1)^{n-i}\binom{i+1}{n-i}2^iC_i, \hskip1.5cm
(\mbox{\cite[A068764]{Sloane}}).
\end{eqnarray*}
Precisely, the number of alternating GNC-trees of $n+1$ points with
exactly $r$ ascents is
\begin{eqnarray*}
\sum_{i+j+k=n}(-1)^{r+k}\binom{i+1}{j}\binom{2i+k}{k}\binom{k}{r-j}2^{i+k}C_i.
\end{eqnarray*}
\end{theorem}
\begin{remark}
When $r=0$, then $j=0$ and there has no alternating GNC-tree of
$n+1$ points with exactly $r=0$ ascents for $n\geq 1$, so we have
\begin{eqnarray*}
\sum_{i=0}^{n}(-1)^{n-i}\binom{n+i}{n-i}C_i=0,   \hskip.5cm (n\geq
1),
\end{eqnarray*}
which is a special case $q=0$ of the Narayana polynomial identity
\cite{MansourSun}
\begin{eqnarray*}
\sum_{i=1}^n\frac{1}{n}\binom{n}{i-1}\binom{n}{i}q^i
&=&\sum_{i=0}^{n}\binom{n+i}{n-i}\frac{1}{i+1}\binom{2i}{i}(q-1)^{n-i}.
\end{eqnarray*}

\end{remark}

When $x=-1,z=1$ in (\ref{eqn 3.4.3}), one has
\begin{eqnarray*}
T_{-1,0,1}^{P}(t)&=& \frac{1+4t-\sqrt{1+8t^2}}{4t}=1-tC(-2t^2) \\
&=&1+\sum_{n\geq 0}(-1)^{n+1}2^nC_nt^{2n+1}.
\end{eqnarray*}
Then we have
\begin{theorem}
The parity of number of alternating GNC-trees of $m$ points
according to the even or odd number of ascents is zero if $m=2n+3$
and $(-1)^{n+1}2^nC_n$ if $m=2n+2$ for $n\geq 0.$
\end{theorem}

\section*{Acknowledgements} The authors are grateful to the
anonymous referees for the helpful suggestions and comments. The
work was supported by The National Science Foundation of China.


\end{document}